\documentclass[12pt]{amsart}
\usepackage{amsmath, color}
\usepackage[all]{xy}
\usepackage{ifpdf}
\ifpdf
 \pdfoutput=1 
 \pdftrue
\fi

	\ifpdf
	\usepackage[pdftex]{graphicx}
	\else
	\usepackage{graphicx}
	\fi

\theoremstyle{plain}
\newtheorem{thm}{Theorem}
\newtheorem{cor}{Corollary}

\newtheorem{lem}{Lemma}
\newtheorem{prop}{Proposition}

\theoremstyle{definition}
\newtheorem{defin}{Definition}
\theoremstyle{remark}
\newtheorem{rem}{Remark}
\theoremstyle{remark}

\def\Q{\mathbb Q}
\def\R{\mathbb R}
\def\slR{\text{SL}(2,\,\R)}

\usepackage{hyperref}

\begin{document}

\ifpdf
	\DeclareGraphicsExtensions{.pdf, .jpg, .tif}
	\else
	\DeclareGraphicsExtensions{.eps, .jpg}
	\fi

\title[Veech surfaces with non-periodic directions]{Veech surfaces with non-periodic directions in the trace field}
\author{ Pierre Arnoux}
\address{Institut de Math\'ematiques de Luminy (UPR 9016),
         163 Avenue de Luminy, case 907,
         13288 Marseille cedex 09,
         France}
\email{arnoux\@ iml.univ-mrs.fr}
\author{Thomas A. Schmidt}
\address{Oregon State University\\ Corvallis, OR 97331}
\email{toms@math.orst.edu}

\date{24 June 2009}
\thanks{The second-named author thanks the Universit\'e Paul C\'ezanne, Marseille and the Max-Planck-Institut f\"ur Mathematik, Bonn for friendly hospitality during the completion of this work.}

 
\begin{abstract}   Veech's original examples of translation surfaces $\mathcal V_q$ enjoying what McMullen has  dubbed ``optimal dynamics'' arise from appropriately gluing sides of two copies of the regular $q$-gon, with $q \ge 3\,$.       We show that every $\mathcal V_q$ whose trace field  is of degree greater than 2  has  non-periodic directions of vanishing SAF-invariant.   (Calta-Smillie have shown that under appropriate normalization,  the set of slopes of directions where this invariant vanishes agrees with the trace field.)     Furthermore,   we give explicit examples of pseudo-Anosov diffeomorphisms whose contracting direction has zero SAF-invariant.   In an appendix, we prove  various elementary results on the containment of trigonometric fields.
\end{abstract}

\maketitle


\section{Introduction}     W. ~Veech's celebrated dichotomy \cite{V} applies in the setting of his initial examples:  Gluing together two copies of a regular $q$-gon by identifying opposite sides by translation results in a translation surface $\mathcal V_{q}$ having exactly two types of directions for linear flow ---completely  periodic or uniquely ergodic.   Indeed, a completely periodic direction $\theta$  enjoys the stronger property that the surface decomposes into a collection of cylinders in the direction, and the moduli of these cylinders are rationally related.  Directly related to this, there is a parabolic $2 \times 2$ real matrix fixing the inverse of the slope, thus $\cot \theta\,$; Veech shows that for these examples, the obvious rotational symmetry and   a parabolic element generate the full (Fuchsian) group of derivatives of the (orientation preserving) affine diffeomorphisms of the surface.    


After Calta \cite{C},  and Calta--Smillie \cite{CS}, one is strongly motivated to normalize by taking a surface in the $\text{GL}(2, \mathbb R)$-orbit of $\mathcal V_{q}$ whose planar presentation has its vertices in the trace field of $\mathcal V_q\,$ --- for then also the set of slopes in which the $\text{SAF}$-invariant vanishes coincides with this field.    Here,  we apply work of A. Leutbecher and others   on cusps of related Fuchsian groups to easily  reach the sobering implication that even surfaces of optimal dynamics have non-parabolic directions  with vanishing $\text{SAF}$-invariant.   

In a certain sense,  the differing planar presentations of the surfaces emphasize different number fields, with the original Veech construction being directly related to $\mathbb Q(\tan \pi/q)$.    Thus the choice of normalized surface requires understanding the relationships of various trigonometric fields.    The literature reveals some minor confusion about these relationships and hence,  in an appendix, we address   the various possibilities for containments of number fields generated over $\mathbb Q$ by trigonometric values.  

\subsection{Results} 

For each of Veech's $\mathcal V_q\,$ we give a $\text{GL}(2, \mathbb R)$-equivalent surface normalized so that its periodic direction field --- in the sense of Calta and Smillie --- equals the trace field.  With this in hand,  we use number theoretic results to find explicit non-parabolic directions with vanishing SAF-invariant.  In particular, we show the following. 
 
\begin{thm}\label{main}   All directions of vanishing SAF-invariant are parabolic on $\mathcal V_q$ if and only if the trace field of $\mathcal V_q$ is at most a quadratic extension of the field of rational numbers.    Thus,   exactly for $q \in \{3, 4, 5, 6, 8, 10, 12\}\,$. 
\end{thm} 

\begin{thm}\label{thmEgs}     The set of indices for which the  surface $\mathcal V_q$ has an affine  pseudo-Anosov diffeomorphism whose contracting direction has vanishing SAF-invariant includes 
$\{7, 9, 14, 18, 20, 24\}\,$.     
\end{thm} 

\begin{rem}\label{remOne}  Our proof  of  Theorem ~\ref{main}  proceeds by first replacing $\mathcal V_q$ by a natural  degree two quotient, clearly identifying normalized surfaces equivalent to the these surfaces 
so as to apply arithmetic results of the school of A. ~Leutbecher on parabolic fixed points of the Hecke groups.   These groups are commensurable, that is conjugate up to finite index, with the Veech groups of the $\mathcal V_q\,$.   This reduces us to a finite number of unresolved  cases; for these, we construct explicit ``special'' hyperbolic matrices   --- and thus affine pseudo-Anosov diffeomorphisms whose contracting direction has vanishing SAF-invariant.         These examples also complete Leutbecher's project.

In fact, the  positive statement  of  Theorem ~\ref{main}   follows from results now already  known: McMullen \cite{Mc2} has shown (using equivalent terminology)
that any Veech surface of quadratic trace field has the property that all its directions of vanishing SAF-invariant are parabolic.       That the $\mathcal V_q$ of higher degree trace field no longer enjoy this property leads us to presume that no Veech surface of cubic or higher trace field  has only parabolic directions have vanishing SAF-invariant.  
  This is another indication that, in genus larger than 2, the situation is much more complicated than the rigid structure  in genus 2 discovered by Calta and McMullen.
\end{rem}

\subsection{Outline} 
In the following section, we briefly recall background.  Further background material appears in the initial piece of  Section ~\ref{secTriSurf}.  That section ends with a presentation of the normalized form of the Veech examples.    In Section ~\ref{secApplyingLeut}, we apply work of Leutbecher and of later authors to reduce to a finite number of cases.    Section ~\ref{secNonParab} rules out the remaining cases by providing explicit examples of special affine pseudo-Anosov maps on Veech surfaces.   Finally, in the Appendix, based largely on the paper \cite{L} of D. H. ~Lehmer, we clarify the relationship between various trigonometric fields.   

\subsection{Supplementary Code}   The second-named author will keep available on his home page a  Mathematica notebook containing the basic computations for the examples for $q= 18, 20, 24$ of Subsection ~\ref{subeven}.   At the time of the submission of this manuscript, the following is a URL for this file ---

\noindent
{\tt http://oregonstate.edu/~schmidtt/ourPapers/Arnoux/AScode.nb}

\subsection{Thanks}   The present note remained in nascent form for several years.   Indeed, the results of Section ~\ref{secApplyingLeut} are referred to in \cite{HS}.    We  have benefited from the prompting of P. Hubert and from discussions with   him,  E. Lanneau and J. Schmurr.

\section{Background} We briefly recall terminology and give more detailed motivation for this work.

\subsection{Translation surfaces}  A {\em translation surface}  is a real surface with an atlas such that, off of a finite number of points, transition functions are translations.  From the Euclidean plane, this punctured surface inherits a  flat  metric,  and this metric extends to the complete translation surface, with (possibly removable)  conical singularities at the punctures.     Due to the transition function being translations,  directions of linear flow on a translation surface are well-defined.

\subsection{Affine Diffeomorphisms and Fuchsian groups}    Post-com\-posing the coordinate function of a chart from the atlas of a translation surface with any element of $\text{GL}^{+}(2, \mathbb R)$ results in a new translation surface.  This action preserves both the underlying topology and the types of the conical singularities.     

Related to this, an {\em affine diffeomorphism} of the translation surface is a homomorphism that restricts to be a diffeomorphism on the punctured flat surface whose derivative is a {\em constant} $2 \times 2 $ real matrix.  W.~ Veech \cite{V}  showed that for any compact translation surface,  the matrices that arise as such derivatives of (orientation- and area-preserving) affine diffeomorphisms form a discrete subgroup of $\text{PSL}(2, \mathbb R)\,$, thus form a Fuchsian group, now referred to as the {\em Veech group} of the surface.    Furthermore,  Veech showed that if this group is a lattice (that is, if it is of finite co-area), then the directions for flow on the translation surface enjoys the aforementioned dichotomy.    Such surfaces are now referred to as {\em Veech surfaces}.

Perhaps  the main reason that W.~Thurston introduced the notion  of affine diffeomorphism is that  it allows visualization of the action of a pseudo-Anosov diffeomorphism, see \cite{T} ---  an affine diffeomorphism is a pseudo-Anosov diffeomorphism exactly
  when its derivative is a hyperbolic matrix (thus of trace greater than 2 in absolute value).    The two real fixed points of this matrix give the expanding and contracting directions of the diffeomorphism.    We refer to these diffeomorphisms as {\em affine pseudo-Anosov} diffeomorphisms. 

\subsection{Sah-Arnoux-Fathi Invariant}   An interval exchange is an interval  map that acts by translating  subintervals.   The Sah-Arnoux-Fathi invariant, or simply {\em SAF-invariant}, of an interval exchange map is $\sum_{i=1}^{n}\, l_i \wedge t_i$ if the   interval
exchange takes intervals of lengths $l_i$ and translates them by the various $t_i\,$.   This invariant defines a homomorphism from the group (under composition) of interval exchange maps to the additive group $\mathbb R \wedge_{\mathbb Q} \mathbb R\,$, hence it vanishes whenever the interval exchange map is completely periodic.  

It can be proved that a minimal exchange map whose  translation lengths generate a rank 2 $\mathbb Z$-module has non-vanishing invariant.  This is in particular the case if the translation lengths belong to a fixed quadratic field; hence any interval exchange map that is quadratic in this sense and has a null invariant is either periodic or reducible, see \cite{A}. Note that this is no longer the case in higher degree:  for any $d>2$, there are minimal exchange maps with lengths in an algebraic field of degree $d$ and a null invariant, \cite{AY}. 

When we consider a flow preserving a transverse measure on a surface, to any tranverse interval we can associate a first return map which is an interval exchange map. It is not difficult to prove that, as long as the interval cuts every orbit of the flow, the invariant does not depend on the choice of the interval, but only on the flow, see \cite{A}. If the flow is minimal, the invariant is the same for all choices of the interval; otherwise, it is the sum of the invariants of the minimal components.

   Kenyon and Smillie \cite{KS} introduced  an invariant for translation surfaces, the $J$-invariant, that takes values in $\mathbb R^2 \wedge_{\mathbb Q} \mathbb R^2\,$.  They showed that appropriate  projections of the $J$-invariant of the surface give  the SAF-invariants associated to the directions $\theta$ and the surface.  If the direction is {\em completely periodic}, in the sense that  every orbit of the  flow is either periodic, or a saddle connection, then the SAF-invariant of this direction vanishes.  In particular, the aforementioned parabolic directions have vanishing SAF-invariant.   

\subsection{Fields: trace, holonomy, periodic direction} 
Calta and Smillie in \cite{CS} show that if a translation surface has at least three directions with vanishing SAF-invariant, then upon appropriately normalizing the surface, the set of slopes of vanishing SAF-invariant forms the union of infinity with a field. They call this the {\em periodic direction field}, and say that the surface is in {\em standard form}.

Each translation surface has a very natural set of directions --- a {\em holonomy} direction for the surface is a  direction for which some closed cycle exists.   Every completely periodic direction is certainly a holonomy direction, but in general the opposite inclusion does not hold: for example, there could be a closed cycle separating two minimal components.    Slightly more subtly,  a {\em holonomy vector} is defined by using the so-called developing map, and gives a flat representative of a cycle.   Kenyon and Smillie \cite{KS} defined the {\em holonomy field} to be the smallest field over which the set of   holonomy vectors is contained in a two-dimensional vector space.    

Gutkin and Judge \cite{GJ} defined the {\em trace field} of a translation surface to be the field extension of the rationals generated by the traces of derivatives of the affine diffeomorphisms of the surface; this is clearly a conjugacy invariant, and does not depend on the choice of a particular translation surface in the $\slR$-orbit. Calta and Smillie \cite{CS} show that a Veech surface has equality of its trace, holonomy and periodic direction fields. 

\subsection{Genus two Veech surfaces: vanishing SAF implies parabolic}  
Results of Calta \cite{C} and of McMullen \cite{Mc}  show that on any  genus two Veech surface every direction of vanishing SAF-invariant is a parabolic direction.     This result uses the fact that the degree of the periodic direction field is bounded by the genus.    Indeed, as already alluded to in Remark ~\ref{remOne}, McMullen \cite{Mc2} shows that if the trace field of a Veech surface is quadratic, then upon appropriately normalizing, the parabolic directions form the full trace field.  

On the other hand, even in genus two,  there are many translation surfaces (with two singularities) with directions that are completely periodic but not parabolic \cite{C}.   Also,  there are examples of genus two surfaces with a direction of vanishing SAF-invariant but for which the surface decomposes into a connected sum of tori, each of which is ``irrationally foliated'' in the given direction \cite{Mc2}. In this case there are two minimal components with nonzero, but opposite, invariant, so that the global invariant vanishes.  Whereas in genus two this  is the only manner to realize a non-completely periodic direction with vanishing SAF-invariant (see the corollary to Theorem 3.7 of  \cite{A2}), in genus three and higher  Arnoux and Yoccoz \cite{AY} already gave examples of affine pseudo-Anosov diffeomorphisms of vanishing SAF-invariant;  see   \cite{HLM} for recent work related to these examples.   For none of these earlier  examples of non-parabolic directions with vanishing SAF-invariant   is the underlying surface a Veech surface.   

\subsection{Special affine pseudo-Anosov diffeomorphisms}   Given any Fuchsian group, the {\em  trace field }is the field extension of $\mathbb Q$ generated by the traces of the elements of the group.    In general, the eigenvalues belong to a quadratic extension of this field, but it can happen that an element of this trace field  is the fixed point of a hyperbolic matrix $M$ in the group.   The standard fixed point formula shows that this is  the case if and only if  $\text{tr}^2(M) -4$ is a square in this trace field. This curious phenomenon has been noted in other contexts,  see in particular \cite{LR}, where they call such $M$ {\em  special hyperbolic}.  
 
     We call an affine pseudo-Anosov diffeomorphism {\em special} if its derivative is a special hyperbolic matrix.   Recall that the larger eigenvalue of this matrix is called the {\em dilatation} of the diffeomorphism.  In general, the dilatation lies in a number field of degree at most twice the genus of the surface.   Usually, this bound is achieved.   However,  a special pseudo-Anosov diffeomorphism has its dilatation in a field of ``unexpectedly'' lower degree.     

  
\section{Normalized triangle surfaces}\label{secTriSurf}  

In the even index $q= 2m$ case, Arnoux-Hubert \cite{AH} note that the interchange of Veech's two regular $q$-gons defines an involution on the surface $\mathcal V_q\,$.   The quotient is a translation surface that is equivalent to $\mathcal M_{q}\,$,  the surface given by identifying by translation opposite sides of the regular $q$-gon gives a surface.    It is easily verified that this quotient surface has the same parabolic directions, and indeed has its Veech group equal to that of $\mathcal V_q\,$.

Independently,  Earle and Gardiner \cite{EG} discussed the $\mathcal M_{q}$ as triangle surfaces.   Let  $a,b,c$ be a relatively prime triple of positive integers and let $q$ be their sum. Let $\mathcal X(a, b, c)$ denote the translation surface arising from the unfolding of a Euclidean  triangle  of angles $a \pi/q, b \pi/q, c \pi/q \,$.  For odd $q$, one has  $\mathcal V_{q} =  \mathcal X(1, 1, q-2)\,$;   Earle-Gardiner show that $\mathcal M_{2m} = \mathcal X(m-1,m-1,  2)$.   

\subsection{Review of construction of triangle surface}

Recall that \cite{KS} discusses the construction of  triangle surfaces.
Let our triangle $T$ be of angles $(\alpha, \beta, \gamma) =  (a \pi/q, b\pi/q,c \pi/q)$ as above, and such that the side opposite $\gamma$ is of length one and agrees with the real line along $[0,1]$.    Thus the vertices of $T$ are $(0, 0), (1,0)$ and, say,  $v\,$.   Using the law of sines,  the vertex $v$  has coordinates  $( t \cos \alpha, t \sin \alpha)$, with $t =  \sin \beta/\sin \gamma\,$.  Calta and Smilie \cite{CS}, proof of Theorem 1.4, prove the following:

\begin{lem} Let  $a,b,c$ be a relatively prime triple of positive integers, and let $q$ be their sum. Then the trace field of $\mathcal X(a,b,c)$ is $\Q(\cos 2\pi/q\,)$.
\end{lem}

Letting $\overline{T}$ denote the triangle whose vertices are the complex conjugates of those of $T$ and $\zeta$ be a primitive $q$th root of unity,   using standard  complex notation,  \cite{KS}  shows that up to appropriate side identifications the surface $(X, \omega)$ can be given as the union $\tilde T = \cup _{j=0}^{q-1}\,  (\zeta^j T)  \;\bigcup\; \cup_{j=0}^{q-1}\,  \zeta^j \overline{T}\,$.   They also show that all vertices of this union are in $\mathbb Q(\zeta)$.

\subsection{Normalizations} 
There are two natural actions of  matrix groups in our setting, each giving a normal form for a translation surface.
\begin{defin}   A translation surface   is {\em projectively normalized} if all slopes of non-vertical homological  directions  lie in the trace field.   A translation surface given as a collection of Euclidean polygons with sides identified by translations is {\em linearly normalized} if all of the vertices have coordinates in the trace field.  
\end{defin}

One easily verifies that a linearly normalized translation surface is projectively normalized.  (For this, see the argument in the proof of  Calta-Smillie Theorem 1.4 in \cite{CS}.)    The theorem itself  can be expressed in the following manner.

\begin{thm}{(Calta-Smillie \cite{CS})}   Any projectively normalized triangle surface has equality of its trace and periodic fields. 
\end{thm}

Taking into account Lemma ~\ref{ksImPrts} of the Appendix,  the aforementioned work of \cite{KS} and \cite{CS} implies the following.  Let 
\[ N := \begin{pmatrix}1&0\\0&1/(\,\sin 2\pi/q\,)\end{pmatrix}\,.\]
\begin{lem}\label{correctNormalization}    Let $a,b,c$ be a relatively prime triple of positive integers and let $q$ be their sum.   Let $\mathcal X(a,b,c)$ denote the translation surface resulting from the standard unfolding process of the triangle of angles $a\pi/q, b\pi/q, c \pi/q$ and of vertices normalized as above.      Let 
\[
\mathcal Y(a,b,c) := \begin{cases} \mathcal X(a,b,c) &\text{if}\; 4 \;\vert\, q\,;\\
\\
N \cdot \mathcal X(a,b,c)&\text{otherwise}\;.
\end{cases}
\]

\noindent
 Then  $\mathcal Y(a,b,c)$  is linearly normalized.  Thus,  the  set of slopes of non-vertical directions of $\mathcal Y(a,b,c)$ whose linear flow has vanishing  SAF-invariant is  $\mathbb Q( \cos 2\pi/q)\,$, the trace field of $\mathcal X(a,b,c)\,$.  
\end{lem} 
J. Schmurr \cite{S} points out that the  authors  of \cite{CS} forgot to apply the above normalization in the proof of their Theorem 1.4 (he gives an alternate proof of their theorem  in his 2008 Oregon State University Ph.D. dissertation).

\section{Applying Leutbecher's characterization of Hecke group cusps}\label{secApplyingLeut}
 
 When the index $q$ is odd, the Veech group  of  $\mathcal V_q$ is a hyperbolic triangle Fuchsian group of signature $(2, q, \infty)$  ---  meaning that (by way of its action on the hyperbolic plane) it uniformizes a once punctured sphere with singularities of torsion index $2$ and $q\,$, respectively.   When the index $q$ is even,  the Veech group is a subgroup of index two in a group of signature $(2, q, \infty)\,$.   
 
 The parabolic directions on a translation surface are in two-to-one correspondence with the parabolic fixed points of the Veech group of the surface (directions $\theta$ and $\theta + \pi$ have identical images).   For brevity, one often refers to the set of the parabolic fixed points of a Fuchsian group, acting on the (extended) reals,  as its {\em cusps}.   Perhaps unfortunately,  it is also traditional to refer to each group orbit of these fixed points as a cusp.   Context usually suffices to clarify this ambiguity.  
  
 We translate our geometric questions into questions about the cusps of a well-studied family of groups that are conjugate to our Veech groups.   In the first subsection, we briefly recall this work.

\subsection{Leutbecher's results}
A standard presentation of the hyperbolic triangle Fuchsian group of signature $(2, q, \infty)\,$, is as the Hecke group 
\[ G_q := \langle \begin{pmatrix} 1 &\lambda\\ 0&1 \end{pmatrix}, \, \begin{pmatrix} 0 &-1\\ 1&0 \end{pmatrix}\,\rangle\,\]
with $\lambda = \lambda_q := 2 \cos \pi/q\,$.

The first of these generators is parabolic, and for its M\"obius action on the extended reals,   has $\infty$ as its fixed point.   From this it follows that every element of the $G_q$-orbit of $\infty$ is a parabolic fixed point.   As the signature indicates, $G_q$ has exactly one cusp:   this orbit is the set of all parabolic fixed points.

Now,   the orbit of infinity   is $G_q\cdot \infty := \{\,  \dfrac{a}{c}\;\vert \,  \begin{pmatrix} a & \star\\ c&\star \end{pmatrix} \in G_q\,\}\,$.   This is obviously a subset of $ \mathbb Q(\lambda) \cup
\{\infty\}\,$, but in fact more can be said.     Here and below,  notation such as $\lambda \, \mathbb Q(\lambda^{2})$ denotes the set obtained upon multiplying all elements of the indicated field by the leading factor. 
\begin{thm}{(Leutbecher \cite{Leu})}\label{thmLeut}  For each  $q\ge 3\, $ 
\[ G_{q} \cdot\infty \subset \lambda \, \mathbb Q(\lambda^{2}) \cup
\{\infty\}\;.\]
Equality holds for $q \in \{3, 4, 5, 6, 8, 10, 12\}\,$.
\end{thm}

\begin{rem} That $G_{q} \cdot\infty$ is a subset of  $\lambda \, \mathbb Q(\lambda^{2}) \cup
\{\infty\}$ is easily shown by induction using the generators displayed above --- an ordered pair $(a, c)$ giving a column of any element of $G_q$ must be such that exactly one element of the pair is in $\mathbb Z[\lambda^2]\,$, and the other is  in $\lambda \mathbb Z[\lambda^2]\,$.  (Note also that only for $q$ even is $ \mathbb Q(\lambda_{q}^{2})$ unequal to   $\mathbb Q(\lambda_q)\,$, and then it is an index two subfield --- see Lemma~\ref{cos2tanSqrd} and Lemma~\ref{lemCos2x} to verify this elementary observation.)

 Leutbecher's list of indices is exactly those for which the  field $\mathbb Q(\lambda_{q}^{2})\,$, thus $\mathbb Q(\cos 2 \pi/q)\,$,  is at most a quadratic extension of $\mathbb Q\,$.  The case of $q=3$ is classical --- indeed,  $G_3 = \text{PSL}(2, \mathbb Z)\,$is  the classical modular group, its parabolic fixed points are certainly infinity and the rationals! ---  $q=4, 6$ are easily derived from this.   For the remaining cases,   Leutbecher's arguments strongly use both the quadratic nature of the fields, and the fact that their rings of integers are principal ideal domains.
\end{rem}

Leutbecher's positive result is complemented by the following.   Let $\lambda \, \mathbb Q(\lambda^2)^{0}$ denote  the subset of  $\lambda \, \mathbb Q(\lambda^{2})$ whose elements may be written as $a \lambda/c$ or $a/ c \lambda\,$ with $a, c \in \mathbb Z[\lambda^2]$ where the pair $a \lambda, c$ (resp. $a, \lambda c$) is relatively prime in $\mathbb Z[\lambda]\,$.     ( Note that this subset can only be proper if $\mathbb Z[\lambda]$ is not a unique factorization domain.)  
\begin{thm}{(Wolfart \cite{W})}  For each  $q\ge 3\, $ 
\[ G_{q} \cdot\infty \subset \lambda \, \mathbb Q(\lambda^{2})^{0} \cup
\{\infty\}\;,\]
Equality holds for at most  $q \in \{3, 4, 5, 6, 8, 9, 10, 12, 18, 20, 24\}\,$.
\end{thm}

\begin{rem} The first part of our remark on Leutbecher's proof,   and the fact that any element of $G_q$ has determinant one,   shows that the containment holds here.   To prove that equality cannot hold for all $q$ not listed,  Wolfart adapts an argument of Borho and Rosenberger \cite{BR}.    Informally,  one replaces $G_q$ acting on $\mathbb Q(\lambda_q)$ by a homomorphic   group that acts on a finite field.   For $q$ not on Wolfart's list,  he uses counting arguments to show that this image group cannot act transitively on the finite field.  However,  this implies intransitivity of the action of $G_q\,$.  
\end{rem}
 
\begin{rem}  The case of odd index  was resolved by \cite{Sei}, who showed that equality does not hold when $q=9\,$.   The examples   from \cite{TetAl} given in Lemma~\ref {lemCandidates} are sufficient for our purposes --- one has that for odd $q$,  the orbit of infinity contains all of $\mathbb Q(\lambda_q)$ exactly for $q=3,5\,$.    

In the even index case, we give explicit (counter)examples for $q=18, 20, 24$ so as to conclude that Leutbecher's list of indices is complete. 
\end{rem}

\subsection{Explicit conjugation} 
Any two triangle groups of the same signature are $\text{PSL}(2, \mathbb R)$-conjugate.   
Arnoux and Hubert \cite{AH} give an explicit conjugation that sends  the Veech group of $\mathcal V_q$ to the Hecke group $G_q$ when $q$ is odd; when $q$ is even, it sends the Veech group  to a subgroup of index 2 of $G_q\,$.      The matrix 
 \[ 
H =  \begin{pmatrix}1& \cos \pi/q\\
\\
                           0 &\sin \pi/q\end{pmatrix}\,  
 \]
acts so as to send the hyperbolic triangle of vertices $\infty, \, - \overline{\zeta_{2q}}, \, \zeta_{2q}$ to the triangle of vertices $\infty,i, i + \cot \pi/q\,$.   From this, one finds that  
 the Veech group of the translation surface  $\mathcal H_q := H^{-1}  \circ \mathcal V_q$ is exactly the Hecke group $G_q\,$, or a subgroup of index two.    

In the remainder of this subsection, we show that  $\mathcal H_q$ is projectively normalized exactly when $q$ is odd.    Hence, one can directly use the Leutbecher school's results to prove Theorem ~\ref{main} in this case.   When $q$ is even, we follow the path from $\mathcal H_q$ to the projectively normalized $\mathcal V_{4k}$ or $N\cdot \mathcal V_{4k + 2}\,$, and find that here also their results are germane.

\subsubsection{Odd index $q$} 
\begin{prop}\label{propAlgSensitiveOdd}  Suppose that  $q$ is odd.     Then 
the Veech surface $\mathcal V_q$ has the property that every direction of vanishing SAF-invariant is parabolic if and only if the parabolic fixed points of $G_q$ are infinity and the field $\mathbb Q(\lambda_q)\,$.   
\end{prop}

This follows,  by a conjugation,  from the following lemma.  
  
\begin{lem}\label{lemHeckeOddNormd}  Suppose that  $q$ is odd.     Then the translation surface $\mathcal H_q$ is projectively normalized.     The cotangents of its parabolic direction fill out all of the projective line of its trace field if and only if the orbit of infinity under $G_q$ is all of this projective line. 
\end{lem}

\begin{proof}  Fix $q$ odd.    Since $\mathcal V_{q} =  \mathcal X(1, 1, q-2)\,$  Lemma~\ref{correctNormalization} gives that  $\mathcal Y(1, 1, q-2)\,$  is linearly normalized.  Thus, the SAF-invariant vanishes for  all directions of flow on this surface with cotangent in $\mathbb P^1 (\,\mathbb Q(2 \cos \pi/q)\,)\,$.       Consider now the matrix product 
\begin{equation}\label{eqMatProd}  
N\, H = \begin{pmatrix}1&\cos \frac{\pi}{q}\\
  \\
  0 &1/( 2 \cos \frac{\pi}{q})\end{pmatrix}\,.
\end{equation} 
Since $q$ is odd,  Lemma ~\ref{lemCos2x} certainly implies that  $\cos \pi/q \in \mathbb Q( \cos 2\pi/q)\,$.  
We deduce that the matrix product is contained in $\text{GL}_{2}(\,\mathbb Q(2 \cos \pi/q)\,)\,$.   Therefore,   acting as a M\"obius transformation, it sends $\mathbb P^1 (\,\mathbb Q(2 \cos \pi/q)\,)\,$   to itself.   It thus follows that also on $\mathcal H_q$  the slopes of directions with vanishing SAF-invariant form $\mathbb P^1 (\,\mathbb Q(2 \cos \pi/q)\,)\,$.     That is,  $\mathcal H_q$ is projectively normalized.   
 
 Now, the parabolic directions of $\mathcal H_q$ have as their cotangents the parabolic fixed points of $G_q\,$, and these form exactly the orbit of infinity under $G_q\,$.    
 \end{proof} 

\begin{rem} In fact, one can show that $\mathcal H_q$ is linearly normalized.
\end{rem}

\subsubsection{Even index $q$} 
 
\begin{prop}\label{propAlgSensitiveEven}  Suppose that  $q$ is even.   The Veech surface $\mathcal V_q$ has the property that every direction of vanishing SAF-invariant is parabolic if and only if the parabolic fixed points of $G_q$ form the set  $\lambda_q \;\mathbb Q(\lambda_{q}^{2}) \cup \{\infty \}\,$.   
\end{prop}

 This follows from the following two lemmas;  note that although the Veech group of  $\mathcal V_{q}$  is of index two in the conjugated copy of $G_q$, it has the same set of cusp values (there is simply a partition of these into two sets, one for each of the cusps of the uniformized hyperbolic surface).

\begin{lem}\label{lemDiv4Cusps}  Suppose that the index $q$ is divisible by 4.    Then 
$H$ sends  $\lambda_q \;\mathbb Q(\lambda_{q}^{2})$ to  $\mathbb Q(2 \cos \pi/q)\,$.
\end{lem}

\begin{proof} For any $q$,  Lemma ~\ref{cos2tanSqrd} implies  that 
\[  2\;  \frac{\cos \pi/q}{\sin  \pi/q} \; \mathbb Q(\cos^2 \pi/q) + \cot \pi/q =   2 \cot \frac{\pi}{q} \;\; \mathbb Q(\tan^2 \pi/q) + \cot \pi/q\;.\]
Since here $q$ is a multiple of 4,   Lemma~\ref{tangents} gives that $\mathbb Q(\tan^2 \pi/q) = \mathbb Q(\tan \pi/q)\,$.    Thus, $H$ sends $ \lambda_q \, \mathbb Q(\lambda_q^{2})$ to 
\[ 
 2 \cot \frac{\pi}{q} \;\; \mathbb Q(\tan  \pi/q) + \cot \pi/q =  \mathbb Q(\tan  \pi/q) \;,
 \] 
 since multiplication by any non-zero element of a field defines a surjection.  Finally,  by Lemma ~\ref{lemCos2x},  since 4 divides $q$,  this field is  $\mathbb Q(\cos  2\pi/q)\,$.  
 \end{proof} 

\begin{lem}\label{lem2Mod4Cusps}  Suppose that the index $q$ is congruent to 2 modulo 4.    Then applying
$N \, H$ sends  $\lambda_q \;\mathbb Q(\lambda_{q}^{2}\,)$ to  $\mathbb Q(2 \cos \pi/q)\,$.
\end{lem}

\begin{proof}  The product $N \, H$  is given in Equation ~\eqref{eqMatProd}, and one easily finds that 
$N \, H$ sends $\lambda_q \;\mathbb Q(\lambda_{q}^{2}\,)$ to 
\[ 2 \cos \frac{\pi}{q}\, \bigg [\, 2 \cos \frac{\pi}{q}\, \mathbb Q( \cos \frac{2 \pi}{q}\,) + \cos \frac{\pi}{q}\,\bigg]\,.\]
Since for any $\theta\,$,  we have $\cos^2 \theta \in \mathbb Q( \cos 2 \theta\,)\,$,  the result easily follows.  
 \end{proof}

\subsection{Proof of Theorem ~\ref{main}}
If $q$ is odd, then Proposition ~\ref{propAlgSensitiveOdd} and Leutbecher's result (Theorem \ref{thmLeut} above) give the first part of Theorem ~\ref{main}.   The proposition combined with Wolfart's and Seibold's results complete the proof in this odd index case. 

When $q$ is even,  an appropriate complement to Wolfart's result is provided by examples in the next section, showing  that Leutbecher had indeed already listed all $q$ for which $G_q \cdot \infty$ contains all of $\lambda \; \mathbb Q(\lambda^2)\,$.    Combining this with 
 Proposition ~\ref{propAlgSensitiveEven} gives Theorem ~\ref{main} in this even index case.


\section{Special affine pseudo-Anosov diffeomorphisms}\label{secNonParab}

Since $\mathbb Q(\lambda_{q}^{2}) = \mathbb Q(\lambda_{q})$ when  $q$ is odd,  
  Leutbecher's theorem here exactly addresses the question of when all elements of the field are parabolic fixed points of $G_q\,$.  For these odd $q$,  any element of the field $\mathbb Q(\lambda_{q})$ that is not in $G_q \cdot \infty$ corresponds to a non-parabolic direction that lies in the trace field of $\mathcal H_q$ and thus has vanishing SAF-invariant.    Explicit examples include  those given in Subsection ~\ref{subOdd}  by Rosen and Towse \cite{RT}, and by Towse {\em et al} \cite{TetAl}.    

Rosen and Towse only treat small index cases,  and for even $q$ Towse {\em et al}  only discuss (the numerous) $G_q$-orbits of $\mathbb Q(\lambda)$, whereas we are interested in the action when restricted to $\lambda \; \mathbb Q(\lambda^{2})\,$.   We give new examples of fixed points of special hyperbolic matrices, and thus also of special affine pseudo-Anosov diffeomorphisms, in Subsection ~\ref{subeven}.  

\subsection{Odd index $q\,$:  non-parabolic points}\label{subOdd}
\begin{lem}\label{lemCandidates}(Rosen-Towse \cite{RT}, Towse {\em et al} \cite{TetAl})   Each of the following elements is not in the orbit of infinity under its respective Hecke group.  
\begin{itemize} 
\item[(i.)] $\;\; q = 7: \;\; \lambda_{7}^2 - 1    \,; $
\smallskip
\item[(ii.)] $\;\;  q = 9: \;\; \pm (2\lambda_{9}^2 +2)\;; \;  \pm (8 \lambda + 8)    \,. $
\end{itemize}
\end{lem} 

\begin{cor}\label{corCandidates}  The above examples are fixed points of special hyperbolic matrices.
\end{cor} 

\begin{proof}  In each of \cite{RT} and \cite{TetAl},  the examples referred to above are shown to have periodic (Rosen) $\lambda_q$-continued fraction expansion.   Now,  any real number of periodic expansion is easily shown to be the fixed point of a matrix in the corresponding Hecke group.  Rosen \cite{R} showed that (1) for each $q\,$,  each real number $x$ has a unique  $\lambda_q$-continued fraction expansion,  and (2) the finite length $\lambda_q$-continued fraction expansions are in one-to-one relationship with the parabolic fixed points of the corresponding Hecke group.   Therefore,  any periodic expansion corresponds a real fixed point of a non-parabolic matrix.   Thus,  a real number of periodic expansion is indeed the fixed point of a hyperbolic element. 
\end{proof}

\begin{rem} Some of the examples of Towse {\em et al.}   were found by also using a variant of the Borho-Rosenberger   arguments.    In specific cases, one can identify  infinite classes of elements in  $\mathbb Q(\lambda_q)$ that cannot be in the orbit of infinity.     These elements can then be examined as to whether they are fixed points of (special) hyperbolic elements (as all of the examples above), or not.    
 
 Note also  \cite{TetAl} conjecture that for every  $q >9$ there are elements of $\mathbb Q(\lambda_q)$ that are fixed by no elements of $G_q\,$.   They also report having searched without success for (in the present terminology) hyperbolic fixed points in $\mathbb Q(\lambda_q)$ for odd $q$ in the range from $11$ to $29\,$.   Based on this, they indicate some doubt as to  the existence of such fixed points for $q>9\,$.
 \end{rem}
\bigskip

\subsection{Even index $q\,$: new examples} \label{subeven}
Towse and his various co-authors focused on the question of expansions of elements of $\mathbb Q(\lambda_q\,)$,  but we are interested in properties of elements of $\lambda_q\; \mathbb Q(\lambda_{q}^{2}\,)$.  Thus, we further examine this  even index setting.  We give explicit examples of special hyperbolic matrices for $q=14, 18, 20, 24\,$.    Note that this also establishes that the corresponding $\mathcal V_q$ have special affine pseudo-Anosov diffeomorphisms.   

We first give an example for the smallest even index for which this is possible.   

\begin{lem}\label{q14Candidate}      Fix $q=14$ and $\lambda = \lambda_{14}\,$.   Let $S = \begin{pmatrix} 1 &\lambda\\ 0&1 \end{pmatrix}\,$ and $T = \begin{pmatrix} 0 &-1\\ 1&0 \end{pmatrix}\,$.    Let $x$ be the attractive  fixed point of  
\[ M := S T S^{-1}TS^{-1} T S T\;.\]
Then $H\cdot x \in \mathbb Q( \cos \pi/7\,)\,$ and thus the conjugate of $M$ by $NH$ is a special hyperbolic matrix.
\end{lem} 

\begin{proof}   For any $\lambda\,$ we find that a matrix of this form gives 
\[ M = \begin{pmatrix} \lambda^4 + \lambda^2 + 1 &-\lambda^3\\ \lambda^3& -\lambda^2 + 1  \end{pmatrix}\,,\]
and is thus certainly hyperbolic, and fixes 
\[ x = \dfrac{ \lambda^2 +2 \pm \sqrt{ \lambda^4+4}}{ 2  \lambda}\,.\]

Now, the minimal polynomial over $\mathbb Q$ of $\lambda_{14}$ is $p(x) := x^6 - 7 x^4 +  14 x^2  -7\,$ and we find that $(x^4-25)^2  \equiv (x^4 + 4) (x^4 - 12)^2 \; \text{mod}\; p(x)\,$.   Thus,  
\[\lambda^4 + 4 \, = \, \big( \dfrac{\lambda^4 - 25}{\lambda^4 - 12}\big)^2\,.\]
We conclude that $ x \in \lambda\,  \mathbb Q(\lambda^2)$ and hence  the result now follows from Lemma~\ref{lem2Mod4Cusps}.   
\end{proof} 

\bigskip 

In briefer form,  we give examples for the three indices for which Wolfart's result does not immediately imply the intransitivity of $G_q$ on $\lambda \mathbb Q(\lambda^2)\,\cup \{\infty\}\,$.  
\begin{lem}\label{q18Candidate}       Let $\lambda = \lambda_{18}\,$, 
 $M  =  STS^4TS^{-1}TS^{-4}T\,,$ 
and let $\delta = \text{tr}^2(M) - 4\,$.    Then 
\[
\delta =  1 + 4 \lambda^4\;\; \text{and} \;\; \sqrt{\delta} = \dfrac{3 \delta + 37}{\delta - 33}\,.  
\]
The matrix $M$ fixes the point
\[\dfrac{7 \lambda^4 - 17 \lambda^2 + 12}{2 \lambda (\lambda^4 - 8)}\,.\]
 \end{lem} 
\begin{proof}   The above is verified by direct calculation.   
\end{proof} 

\bigskip

\begin{lem}\label{q20Candidate}    
    Let $\lambda = \lambda_{20}\,$, and 
\[
M  :=  S^4TS^{-4}T(S^{-1}T)^{3}\,,
\]
and let 
\[\delta = \dfrac{\text{tr}^2(M) - 4}{\lambda^2 (-2 + \lambda^2)^2}\,.\]
  Then 
\[
\delta = 1 + 4 \lambda^4-4 \lambda^6 + \lambda^8
\]
and 
\[
\sqrt{\delta} =   \dfrac{\delta^2 - 29  \delta + 41}{14 \delta -16}\,,
\]
and $M$ fixes the point
\[x := \dfrac{-19 + 202\lambda^2 - 273 \lambda^4 + 65 \lambda^6}{\lambda\, (8   + 
   40 \lambda^2 - 68 \lambda^4 + 18\lambda^6\,)}\,.\]
   Furthermore,    $H\cdot x \in \mathbb Q( \cos \pi/10\,)\,$ and thus the conjugate of $M$ by $H$ is a special hyperbolic matrix.

\end{lem} 

\begin{proof}  Again, verification of the claimed equalities is straightforward.   Lemma ~\ref{lemDiv4Cusps} finishes the proof.
\end{proof} 

\begin{rem} 
Whereas the special hyperbolics as above were found by searching amongst short words in our chosen generators,  the next example was found by searching for elements of $\lambda_{24} \mathbb Q(\lambda_{24}^{2})$ with periodic Rosen continued fraction expansions.   Indeed, $x =(2 \lambda^2 + 13)/(\lambda + \lambda^3)$ is in the $G_{24}$-orbit of the fixed point of $M$ below.  Thus, of course it too is fixed by a special hyperbolic matrix, and hence has eventually periodic expansion.  However, as the preperiod is quite long,  we prefer to give $M\,$ corresponding directly to the period. 
\end{rem}

\begin{lem}\label{q24Candidate}        Fix $q=24$ and $\lambda = \lambda_{24}\,$.     Let
\[
W = (TS)^4 T S^6 (TS^{-1})^2 T S^{-5} (TS)^2
\]
and 
\[
W' = (TS^{-1})^4 T S^{-6} (TS)^2 T S^{5} (TS^{-1})^2.  
\]

  Then 
  \[ M := W\, W'\;\]
   is a special hyperbolic matrix.
\end{lem} 


\begin{proof} The minimal polynomial of $\lambda$ is $m(x) =  1 - 16 \lambda^2 + 20 \lambda^4 - 8 \lambda^6 + \lambda^8\,$.   Taking polynomial remainders with respect to the modulos $m(x)\,$,  one thus finds that  $M = \begin{pmatrix}a&b\\c & d\end{pmatrix}\,$ with

\[
\begin{aligned}
a &= -351899 + 5540990 \lambda^2 - 5630280 \lambda^4 + 
   1384810 \lambda^6\,,\\
   b &=  
  267120 \lambda - 4206030 \lambda^3 + 4273780 \lambda^5 - 
   1051090 \lambda^7\,,\\
   c &= 216120 \lambda - 3402520 \lambda^3 + 
   3454760 \lambda^5 - 848800 \lambda^7\,,\\
   d &= 
   -643899 + 
   10138610 \lambda^2 - 10297720 \lambda^4 + 2530790 \lambda^6\,.
   \end{aligned}
   \]
One directly verifies that $M$ has trace (much!) larger than 2, and  fixes 
\[ \dfrac{-38424 \lambda + 599662 \lambda^3 - 527439 \lambda^5 + 95988 \lambda^7}{-23909 + 
   374995 \lambda^2 - 338345 \lambda^4 + 62189 \lambda^6}\,,
  \] 
and thus $M$ is special. 
\end{proof}

\begin{rem}  Our examples for $q \in \{18, 20, 24\}$ complement Wolfart's (and Seibold's) results and thus complete the first part of the Leutbecher program --- it is exactly for the list in Theorem ~\ref{thmLeut} that the orbit of infinity under $G_q$ contains all of $\lambda\, \mathbb Q(\lambda^2)\,$.      A second, yet more difficult,  goal is to arithmetically characterize the cusp set for any $q\,$.    Note that this second goal is also of great interest in our setting --- one would like to precisely locate the set of parabolic directions amongst all those of vanishing SAF-invariant.  It is also completely motivated in the initial setting:  the Rosen $\lambda_q$-continued fractions define a numeration system;   one of the initial problems in the study of any enumeration system is to determine the set enumerated by finite length expansions (for these continued fractions, this is indeed the orbit of infinity, as Rosen showed).    
 \end{rem}
\bigskip



\section{Class number restriction}  Leutbecher and his school were well aware that 
the orbit of infinity under $G_q$ can contain all of $\lambda_q\, \mathbb Q(\lambda_{q}^{2})$ only if the ring of integers of $\mathbb Q(\lambda_{q}^{2})$ is principal.   
With this as hindsight, one can easily recognize that $\mathcal V_q$ with non-parabolic directions with vanishing SAF-invariant  {\em must} arise.  For ease of discussion, let us consider the odd index case.    Then $G_q$   is a subgroup (of infinite index when $q>3$) of the Hilbert modular group of  the totally real field  $\mathbb K := \mathbb Q(\zeta_{q} + \zeta_{q}^{-1}\,)\,$.    This Hilbert modular group acts on the field $\mathbb K$ union infinity, the number of orbits equals the class number of the ring of integers of the field, $h^{+}_{q}\,$ --- for this classical result, see say \cite{vdG}.    Of course, no subgroup could act so as to give fewer orbits.   However,  $G_q$  has exactly one cusp, corresponding to one orbit of field elements under this group.     Thus, as soon as the class number of the field exceeds two there are certainly elements of the field that are not in the cuspidal orbit of this group.    Using the normalization as above, one deduces that then $\mathcal V_q$ has non-parabolic directions with vanishing SAF-invariant.   

Now, although the class numbers of fields of the form $\mathbb Q(\zeta_q + \zeta_{q}^{-1}\,)$ are ``notoriously difficult to compute'' \cite{Sch},    one does know that for $q=163$ the above argument must give an example:  van der Linden \cite{vdL} shows that both   $h^{+}_{163} \ge 4\,$, with equality assuming the Generalized Riemann Hypothesis;  and for all primes  $q< 163$,  $h^{+}_{q}  = 1\,$, under the same assumption.     




\section*{Appendix:   Containments amongst Trigonometric Number Fields}\label{App}

There is some minor confusion in the literature about containments amongst fields defined over $\mathbb Q$ by adjoining values of sines and cosines.  For example, even the redoubtable computational number theorist \newline D. ~H.~Lehmer published inaccurate  results on this matter \cite{L}.   (That there is a flaw in his work is highlighted by a displayed table  indicating that there is no integer $m$ such that the algebraic number $\sin 2 \pi/m$ is  of degree  three:   both $m = 28, 36\,$are of this type.)    We address this here, of course we acknowledge that these elementary results are well-known.  

Even in the case of a rational multiple of $\pi$,  one finds that various configurations of containments  of trigonometric fields arise.   In this appendix, we prove the following result. 

\begin{thm}\label{thmTrigFldCont}  Assume that $q>2$ is an integer and that $q \neq 4$.  Then we have the following 
\begin{enumerate}
\item  If $8 \vert q$ then 
\[   
\mathbb Q(\cos \dfrac{2 \pi}{q}) = \mathbb Q(\sin \dfrac{2 \pi}{q})\;;
\]
\item  If $q \equiv 4 \bmod 8$ then 
\[  
\mathbb Q(\sin \dfrac{2 \pi}{q}) \; \text{is a subfield of index two in }\; \mathbb Q(\cos \dfrac{2 \pi}{q})\;; 
\]
\item  If $4 \not{\vert} \, q$ then 
\[  
\mathbb Q(\cos \dfrac{2 \pi}{q}) \; \text{is a subfield of index two in }\; \mathbb Q(\sin \dfrac{2 \pi}{q})\;. 
\]
\end{enumerate}

\end{thm}

%

\subsection{Degrees using right triangle relation}
The following is a correction of Lehmer's \cite{L} results, we sketch the proof in a series of lemmas.     
\begin{prop}\label{propCompDeg}   Assume that $q>2$ is an integer and that $q \neq 4$.   Then

\[  
 [\, \mathbb Q(\sin \dfrac{2 \pi}{q})\,:\,\mathbb Q\,]   = \begin{cases}  2 \cdot  [\, \mathbb Q(\cos \dfrac{2 \pi}{q})\,:\,\mathbb Q\,] & \text{if}\; 4 \not{\vert} \, q \,;\\
 \\
                                                              \dfrac{1}{2}\cdot  [\, \mathbb Q(\cos \dfrac{2 \pi}{q})\,:\,\mathbb Q\,]  & \text{if}\; q \equiv 4 \bmod 8\,;\\
                                                              \\
                                                               [\, \mathbb Q(\cos \dfrac{2 \pi}{q})\,:\,\mathbb Q\,]  & \text{if}\; 8 \; \text{divides}\;  q\,.
                                   \end{cases}                           
\]

\end{prop}

\begin{lem} \label{lemCosDeg} With $q>2$,     
\[  
[\,  \mathbb Q(\cos \dfrac{2 \pi}{q}) :  \mathbb Q\,] = \phi(q)/2\,.
\]
\end{lem} 
\begin{proof}    Recall that  for any $n \in \mathbb N$ letting  $\zeta_n = e^{2 \pi i/n}\,$   gives $[\mathbb Q(\zeta_{n}): \mathbb Q)\,] = \phi(n)\,$, where as usual $\phi(n)$ denotes the Euler totient function of $n$.     But,  $\zeta_q + \zeta_{q}^{-1} = 2 \cos 2 \pi/q$ and also $\mathbb Q(\zeta_q + \zeta_{q}^{-1})$ is the fixed field of complex conjugation on $\mathbb Q(\zeta_q)\,$.
\end{proof}

\begin{lem} \label{lemRelPrime} With $q>2$,     if $k\in \mathbb N$ is  relatively prime to $q\,$, then
\[  
 \mathbb Q(\cos \dfrac{2 \pi k}{q}) =  \mathbb Q(\cos \dfrac{2 \pi}{q})\,.
\]
\end{lem} 
\begin{proof}     The map $\zeta_q \mapsto \zeta_{q}^k$  defines a Galois automorphism of the cyclotomic field $\mathbb Q(\zeta_q)\,$ that descends to $ \mathbb Q(\cos \dfrac{2 \pi}{q})\,$.
\end{proof} 
 
The following is an immediate implication of the fact that $\sin x$ and $\cos x$ are complementary trigonometric functions. 
\begin{lem} \label{lemSin} With $q>2$,    
\[  
 \mathbb Q(\sin \dfrac{2 \pi}{q}) =  \mathbb Q(\cos \dfrac{2 \pi (4-q)}{4q}\, )\,.
\]
\end{lem} 

The previous three lemmas imply the following.
\begin{lem} \label{lemDeg} With $q>2$,    
\[
[\,  \mathbb Q(\sin \dfrac{2 \pi}{q}) \,:\, \mathbb Q\,] = \dfrac{1}{2}\, \phi(K\,)\,,
\]
where

\[ K = \dfrac{4q}{\gcd(\,4q, 4-q\,)}\,.\]
\end{lem} 

 Lehmer failed to notice the distinction of the final three cases in the following result.  
\begin{lem} \label{lemValue} With $q>2$,    
\[
 \gcd(\,4q, 4-q\,) = \, \begin{cases}   1 & \text{if}\;\; q \; \text{is odd}\,;\\
                                                              2 & \text{if}\;\; q \equiv 2 \bmod 4\,;\\
                                                              4 & \text{if}\;\; 8 \; \text{divides}\;  q\,;\\
                                                             8 & \text{if}\; \;q \equiv 12 \bmod 16\,;\\
                                                             16 & \text{if}\; \; q \equiv 4 \bmod 16\,.
                                   \end{cases}                           
\]
\end{lem} 

\begin{proof}   The cases  odd $q$  and $q \equiv 2 \bmod 4$ are completely trivial.  If $4$ divides $q$, we write $q = 4n$, and have that $ \gcd(\,4q, 4-q\,) = 4  \gcd(\,4n, 1-n\,)\,$.  Now,  

\[ \gcd(\,4n, 1-n\,)  = \, \begin{cases}   1 & \text{if}\; n \; \text{is even}\,;\\
                                                              2 & \text{if}\; n \equiv 3 \bmod 4\,;\\
                                                            4& \text{if}\; n \equiv 1 \bmod 4\,.
                                   \end{cases}                           
\]
Our result follows.
\end{proof} 
The previous four lemmas combine to give the proof of Proposition ~\ref{propCompDeg}.    Indeed, $K$ takes the values $4q$, $2q$, $q$, $q/2$ and $q/4$ respectively in the five cases for $q$ of the previous lemma.   Thus $\phi(K)$ takes the value:   $2 \phi(q)$ when $q$ is either odd,  or is  $2 \bmod 4$; $\phi(q)$ when $8$ divides $q\,$;  and $\phi(q/2) = \phi(q/4)= \phi(q)/2$ in the final two cases, as there $q \equiv 4 \bmod 8\,$.

\subsection{Containments using Pythagorean relations}  
We finish the proof of Theorem ~\ref{thmTrigFldCont} by the following yet more elementary observations.

\begin{lem}\label{cos2tanSqrd} For  $x \in \mathbb R\setminus \{\pi/2 +\pi
\mathbb Z\}$,
\[{\mathbb Q}( \cos\, 2 x\,) = {\mathbb Q}(\tan^2\, x)= {\mathbb Q}(\cos^2\, x)= {\mathbb Q}(\sin^2\, x)\; .\]
\end{lem}

\begin{proof}    We use the identities $\tan^2 x + 1 = \cos^2 x\,$, $ \cos 2 x =  2 \cos^2 x - 1\,$ and $\sin^2 + \cos^2 = 1\,$.
\end{proof}

\begin{figure}
\begin{center}
\includegraphics{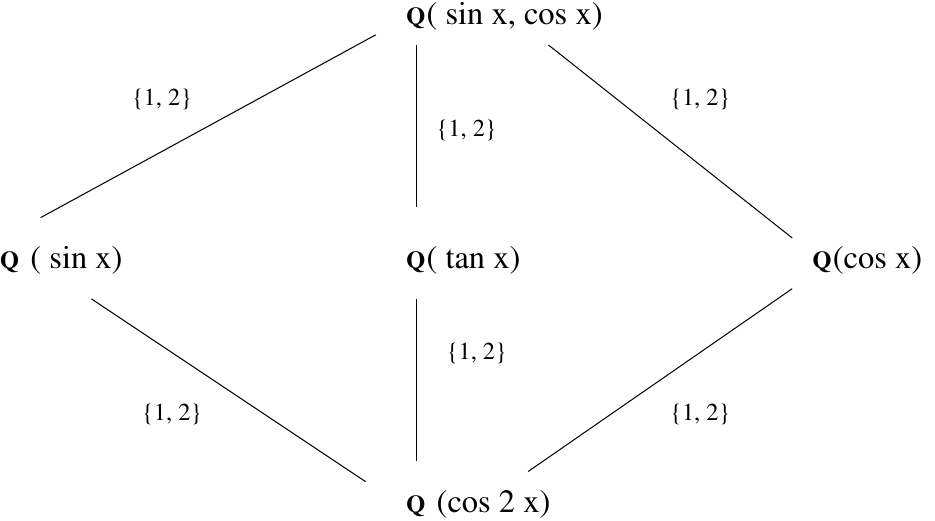}
\caption{Basic field diagram for trigonometric fields. Each intermediate extension degree is either one or two.}
\label{trigFieldsFig}
\end{center}
\end{figure}

\begin{lem}\label{lemTanTwos} For  $x \in \mathbb R\setminus \{\pi/2 +\pi
\mathbb Z\}$,  Figure ~\ref{trigFieldsFig} is a valid diagram of extension fields.
\end{lem}

\begin{proof}      For any field $k$ and any $x$, we have that $[\, k(x)\,:\, k(x^2)\,] \in \{1, 2\}\,$.   Thus, Lemma ~\ref{cos2tanSqrd} shows that the bottom row of field extensions are correctly indicated.  Now, 
clearly each of  $\sin \, x,\, \cos\,  x$ and $\tan \,x$ is contained in ${\mathbb Q}(\sin \, x,\, \cos\,  x\,)\,$.     Using $\sin^2 + \cos^2 = 1\,$,   and its quotients by   $\sin^2 x$ and $\cos^2 x\,$ respectively, shows that also the top row of field extensions is valid. 
\end{proof}

\subsection{Applications}

We use the following easy results in the body of this paper.
\begin{lem}\label{lemCos2x}  With $q>2\,$ 
if $q$ is odd then $\mathbb Q( \cos 2\pi/q) = \mathbb Q(\cos  \pi/q)\,$; if $q$ is even then 
$\mathbb Q(\cos 2 \pi/q)$ is  a subfield of index two in $\mathbb Q(\cos  \pi/q)\,$.
\end{lem}
\begin{proof}  This is an immediate consequence of Lemma ~\ref{lemCosDeg} and Lemma ~\ref{cos2tanSqrd}.
\end{proof}

\begin{lem}\label{tangents} If $q$ is divisible by 4, then 
\[ {\mathbb Q}(\tan^2\, \pi/q)= {\mathbb Q}(\tan\, \pi/q)\; .\]
\end{lem}

\begin{proof}     Lemmas~\ref{cos2tanSqrd} and  ~\ref{lemTanTwos} show that the containment 
$ {\mathbb Q}(\tan\, x) \supset  {\mathbb Q}(\cos\, 2 x) =  {\mathbb Q}(\tan^2\, x)\,$ always holds.   They also show that any time the containment
${\mathbb Q}(\cos\, 2 x) \supset {\mathbb Q}(\sin\, 2 x)\,$ holds, then we must also have the equality $ {\mathbb Q}(\tan\, x)   =  {\mathbb Q}(\tan^2\, x)\,$.   

 But by Theorem ~\ref{thmTrigFldCont}, we have that ${\mathbb Q}(\cos\, 2 x) \supset {\mathbb Q}(\sin\, 2 x)\,$ whenever $x = \pi/q$ and 4 divides $q\,$.
\end{proof}

We  also correct a  typographic error in a statement of \cite{KS}.
\begin{lem}\label{ksImPrts}   Let $\alpha \in \mathbb Q(\zeta_q)$ then 
\begin{enumerate}
\item  If $4 \vert q$ then the imaginary part of $\alpha$ lies in $\mathbb Q(\cos \dfrac{2 \pi}{q})\,$;
\item Otherwise,   the imaginary part of $\alpha/ \sin \frac{2 \pi}{q}$ lies in $\mathbb Q(\cos \dfrac{2 \pi}{q})\,$.
\end{enumerate}
\end{lem} 
\begin{proof} Recall that $\{ \zeta_{q}^{j}\,|\, 1 \le j < q, \gcd(j,q) = 1\}$ gives a $\mathbb Q$-vector space basis of $\mathbb Q(\zeta_q)\,$.     Thus, give  $\alpha \in \mathbb Q(\zeta_q)\,$, there are $a_j \in \mathbb Q$ such that $\text{Im}(\alpha) = \sum\, a_j \, \text{Im}\, \zeta_{q}^{j} = \sum\, a_j \sin 2 \pi j/q\,$.   Now,  in our arguments above, we can replace $\zeta_q$ by any of these $\zeta_{q}^{j}\,$.   Thus,  when $4$ divides $q$,   the first two parts of Theorem ~\ref{thmTrigFldCont} imply that $\mathbb Q(\sin {2 j  \pi}/{q}) \subset \mathbb Q(\cos {2 \pi}/{q}) \,$, and the first part of this lemma follows.\\ 

Recall that 
\[\dfrac{\sin (n+1) \theta}{\sin \theta} = U_n( \cos \theta)\;,\]
where $U_n(x)$ is a Chebyshev polynomial of the second type.  Since the  $U_n(x)$ are polynomials with integer coefficients, the second part of the lemma also follows.
\end{proof}

\end{document}